\documentclass[12pt]{article}
\usepackage{geometry}
\usepackage{amssymb}
\usepackage{color}

\usepackage{amsmath}

\newif\ifleft
\lefttrue 

\newif\ifpre
\pretrue 

\newif\ifhide
\hidefalse

\newtheorem{lemma} {Lemma}
\newtheorem{theorem} {Theorem}

\newif \ifshowup
\showupfalse 

\begin{document}

\title{On mean-field \(GI/GI/1\) queueing model: 
existence and uniqueness
\thanks{
This study has been funded by the Russian Academic Excellence Project '5-100' (section 2.2) and by the RFBR grant \mbox{17-01-00633$\_$a (section 2.1)}.
}
}

\author{A.Yu. Veretennikov\footnote{
School of Mathematics, University of Leeds, Leeds, LS2 9JT, UK; 
email: a.veretennikov@leeds.ac.uk
              \&
              National Research University Higher School of Economics, Moscow, Russian Federation, 
              \&
              Institute for Information Transmission Problems, Moscow, Russian Federation
}         
}

\maketitle

\begin{abstract}
A mean-field extension of the queueing system \(GI/GI/1\) is considered. The process is constructed as a Markov solution of a martingale problem. Uniqueness in distribution is also established under a bit different sets of  assumptions on intensities. 

\noindent
keywords: GI/GI/; mean-field; existence; weak uniqueness; 
Skorokhod lemma
MSC: 60-02; 60K25; 90B22
\end{abstract}

\section{Introduction}\label{intro}
Mean-field approach  in the theory of queueing systems allows to take into consideration large interacting ensembles of queues by using the idea of replacing these interactions by a suitable ``mean field''. This approach showed fruitful in systems with countable and more general state spaces, see, for example, \cite{Aghajani}, \cite{Baccelli}, \cite{Ryb2}, \cite{borovkov}, \cite{Dawson},  \cite{Ryb1}, and the  references therein. However, 
to the best of the author's knowledge, so far there was no general method of constructing mean-field extensions of a basic queueing model such as $GI/GI/1$ -- or, more precisely, $GI/GI/1/\infty$ -- in the literature. In this work we propose such a method   under certain restrictions on intensities of arrivals and service, which intensities may both depend on the state as well as on the marginal distribution of the process. 
This kind of dependence is natural in the study of limits of so-called multi-agent systems, or in other words, of a large number of weakly interacting queues (cf. again the references \cite{Aghajani}, \cite{Baccelli}, \cite{Ryb2}, \cite{borovkov}, \cite{Dawson},  \cite{Ryb1}); in this paper we do not study such a setting  because it seems reasonable to separate the problem of convergence from the problem of existence of the limiting process.  
Existence and weak  uniqueness are discussed on the basis of compactness of measures, Skorokhod's unique probability space Lemma, total variation metric and a Skorokhod--Girsanov's density of measures theorem for jump processes. The basis for the study in the sections \ref{sec:2.1} and \ref{sec:2.2} is a technique similar to the one developed in the preprint on McKean-Vlasov stochastic equations \cite{MV}. Note that in some earlier papers and monographs intensities of transitions in queueing systems were assumed to depend only on the number of customers in the system. This means that the (conditional) distributions  of the service times as well as the arrival times are exponential.  In some situations this is not realistic. In particular, it does not allow heavy tails. Hence, motivation of our extension to a more general dependence is simple: it should relax the assumption of exponential arrival and service times. In the Theorem \ref{thm1} below heavy tails of (conditional) service time distributions are allowed. In the Theorem \ref{thm2} currently heavy tails are not possible, although, both service and arrival times still can be more general than exponential.  Nevertheless, the author's belief is that this is a technical matter to include the case of heavy tails in the conditions for uniqueness, too, which problem will be hopefully resolved in the near future. 
It is likely that the established results may be useful in the area of mathematical theory of reliability which is notably known to have the same basic formulae as queueing theory, see \cite{GBS}. 

The paper consists of Introduction, Main section and t
So, we can rigorously apply  the principle of ``complete probabilities'' (or, better ``complete expectation'') taking a summation as earlier in the intuitive version:wo Appendices. The Main section consists of two subsections related to the two topics shown in the title,  with one theorem in each and with the proof of this theorem. The Appendix 1 contains the statement of Skorokhod's Lemma about an equivalence of weak convergence of a sequence of processes to a convergence in probability of processes with the same distributions on a unique probability space, included for the reader's convenience. The Appendix 2 offers a strict version of a non-rigorous calculus in the middle of the proof of the Theorem~\ref{thm1}.

\section{Main section}\label{sec:1}
The state space of the process under consideration is the union 
\[
{\cal X}:= 
(0,x) \cup \bigcup_{k=1}^{\infty} (k,x,y), \quad x,y\ge 0. 
\]
The meaning of $k$ here is the number of ``customers'' in the system; the value $x$ stands for the elapsed time from the last arrival, while \(y\) signifies the elapsed time of the current service. There is only one server which works without breaks (if there is at least one customer in the system) and it is always in a working state. All newly arrived customers stand in a queue of the infinite capacity, and for simplicity only we assume the FIFO discipline of service (First In First Out). It is assumed that at any time $t$ at any state $X=(k,x,y)$ (or $X=(0,x)$ for $k=0$) there are {\em intensities} of service $\Lambda^-[t,X_t,\mu_t]$ and arrivals $\Lambda^+[t,X_t,\mu_t]$, where $\mu_t$ is the distribution of the random variable $X_t$ itself. Note that occasionally we will be using notation \((0,x,y)\) where \(y\) is a ``false'' variable, i.e., we identify all such triples with any \(y\) with a couple \((0,x)\). It will be sometimes convenient to denote $k=k(X), x=x(X), y=y(X)$ for $X=(k,x,y)$. For technical reasons it is convenient to define the distance between two states $X=(k,x,y)$ and $X'=(k',x',y')$ as  
\[
\rho(X,Y):= |k-k'| + |x-x'| + |y-y'|.
\]
The process is piecewise--linear  Markov (PLMP, see \cite{GK}), which simply means that between any two subsequent jumps the continuous components -- \((x,y)\) if \(k>0\), or just \(x\) if \(k=0\) -- grow linearly with rate 1, while the discrete component \(n\) remains unchanged.

~

{\bf The assumptions}: 
\begin{itemize}

\item[(A1)]
There are Borel measurable, non-negative and  bounded functions $\lambda^+(t,X,Y)$ and $\lambda^-(t,X,Y)$.

\item[(A2)]
\[
\Lambda^\pm[t,X,\mu] = \int \lambda^\pm(t,X,Y)\mu(dY)
\]
(NB: Automatically, both $\Lambda^\pm$ are Borel functions of $(t,X)$.)

\item[(A3)]
The functions $\lambda^\pm(t,X,Y)$ are continuous in all variables. 

\item[(A4)]
The functions $\lambda^\pm(t,X,Y)$ are uniformly bounded away from zero except for $\lambda^-(t,(0,x),(n,y,y'))=0$, for any $x,y,y'\ge 0$  (no jump down from any state with zero customers).

~

Let us emphasize that neither Lipschitz nor any other regularity of the intensities $\lambda^\pm$ is assumed, except for continuity in (A3). Probably continuity in $t$ may be relaxed. Note that functions of time and two state variables will be denoted with round brackets, e.g., as $\lambda^\pm(t,X,Y)$, while functions where the third variable is a measure will be written with square brackets like  $\Lambda^\pm[t,X,\mu]$.

\end{itemize}

\noindent
In particular, intensities $\Lambda^\pm$  {\em may} include additional (non-negative) terms not depending on the measure, say, $\lambda^\pm_0(t,X)$; this may  be helpful so as to justify the assumption (A4), as the terms $\lambda^\pm_0(t,X)$ can be reasonably assumed  uniformly bounded away from zero. Emphasize that (A4) will only be used in the Theorem \ref{thm2}, and as was mentioned earlier, there is a plausible hypothesis that this assumption even in this theorem may be relaxed; 
however, we postpone this issue till further investigations.

~

For \(X\in {\cal X}\) let us denote
\begin{align*}
& X^+ := (k+1,0,y), \qquad \mbox{for}\; X = (k,x,y), \; k\ge 0, 
 \\
& X^- := (k-1,x,0), \qquad \mbox{for}\; X = (k,x,y), \; k\ge 1, 
 \\
& X+\delta := (k,x+\delta, y+\delta), \quad \mbox{for}\; X = (k,x,y), \; k\ge 0,
\end{align*}
Naturally, \(X^-\) is not defined for \(X=(0,x)\). 

\subsection{Existence}\label{sec:2.1}
The initial value $X_0$ of the process may be distributed, which distribution is denoted by $\mu_0$ (in particular, $\mu_0$ may be a delta-measure concentrated at one point). 

\begin{theorem}\label{thm1}
Let the assumptions (A1)--(A3) be satisfied. Then for any initial distribution \(\mu_0\) on \({\cal X}\), on some probability space there exists a Markov process \((X_t, \, t\ge 0)\) with marginal distributions \(\mu_t\) and intensities $\Lambda[t,X_t,\mu_t], \, H[t,X_t,\mu_t]$; in other words, such that for any bounded continuous function \(g(X)\) with bounded continuous derivatives in \((x,y)\), the expression
\begin{equation}\label{df1}
M_t := g(X_t) - g(X_0) - \int_0^t L[s,X_s,\mu_s]g(X_s)\,ds
\end{equation}
is a martingale, where for \(X=(k,x,y)\), \(X'=(k',x',y')\), \(k\ge 0\), \(t\ge 0\), 
\begin{align*}\label{gen}
&\displaystyle L[t,X',\mu] g(X)  := \Lambda^+[t,X',\mu](g(X^+) - g(X)) 
 \nonumber \\ \nonumber \\ 
&\displaystyle + 1(n >0)\Lambda^-[t,X',\mu](g(X^-) - g(X)) 
 \\\\ 
&\displaystyle + \frac{\partial}{\partial x}g(k,x,y) 
+ 1(k>0)\frac{\partial}{\partial y}g(k,x,y). 
\end{align*}
Moreover, for any given measure-valued function \((\mu_s, \, s\ge 0)\) in \(L[s,X_s,\mu_s]\), the {\em martingale problem}  (see \cite{EthierKurtz}) \eqref{df1} has a (weakly) unique solution. 
\end{theorem}
The processes \((X^{}_t,\, t\ge 0)\), or later in the proof of the Theorem \ref{thm1} \((X^{n}_t,\, t\ge 0)\) for \(n\ge 1\) being constructed, let us introduce on some probability space {\em independent}  equivalent processes \((\xi^{}_t, \, t\ge 0)\), or,  respectively, \((\xi^{n}_t, \, t\ge 0)\); let  \(\mathbb E'\) stand in all cases  for the integration with respect to the {\em third variable}, e.g., 
\[
\mathbb E' \lambda^\pm(t,X_t^{},\xi^{}_t) := 
\int \lambda^\pm(t,X_t^{},Y)\mu_t^{}(dY),
\]
or
\[
\mathbb E' \lambda^\pm(t,X_t^{n},\xi^{n}_t) := 
\int \lambda^\pm(t,X_t^{n},Y)\mu_t^{n}(dY),
\]
where $\mu^n_t$ is the distribution of $X^n_t$; this will be repeated in the proof of the Theorem. 

~

Recall that 
\[
\Lambda^\pm[t,X',\mu] =  \int \lambda^\pm(t,X',y)\,\mu(dy) 
= \mathbb E' \lambda^\pm(t,X',\xi), 
\]
where $\xi$ has distribution $\mu$. So, the operator $L$ can be also presented in the form, 
\[
L[t,X',\mu] g(X) = \mathbb E' L(t,X',\xi) g(X), 
\]
for $X=(k,x,y)$, where
\begin{align*}\label{gen}
&\displaystyle L(t,X',y) g(X)  := \lambda^+(t,X',y)(g(X^+) - g(X)) 
 \nonumber \\ \nonumber \\ 
&\displaystyle + 1(k >0)\lambda^-(t,X',y)(g(X^-) - g(X)) 
 \\\\ 
&\displaystyle + \frac{\partial}{\partial x}g(k,x,y) 
+ 1(k>0)\frac{\partial}{\partial y}g(k,x,y). 
\end{align*}

~

Note that the (possibly extended) {\em generator} of the Markov process $X_t$ (cf., for example, \cite{EthierKurtz}) is, of course, $L[t,X,\mu_t]$; different variables $X$ and $X'$ in the definition above are needed only for the convenience of the proof.  
Equivalently, {\em Dynkin's identity} holds true for any function \(g(X)\) from the same class, 
\begin{equation}\label{df2}
\mathbb E_{0,X_0} g(X_t) =  g(X_0) +  \mathbb E_{0,X_0} \int_0^t L[s,X_s,\mu_s]g(X_s)\,ds. 
\end{equation}
Moreover, equivalently, for any \(0\le t_1 < t_2 \ldots < t_{m+1}\), and for any Borel bounded functions \(\phi_k(X), \, X\in {\cal X}\), 
\begin{equation}\label{df22}
\mathbb E_{0,X_0} \left(g(X_{t_{m+1}}) -  g(X_{t_m}) -   \int\limits_{t_m}^{t_{m+1}} L[s,X_s,\mu_s]g(X_s)\,ds\right)
\prod_{k=1}^{m} \phi_k(X_{t_k}) = 0. 
\end{equation}
Also note that for the validity of the equation (\ref{df22}) for any Borel bounded continuous functions \(\phi_k(X), \, X\in {\cal X}\) it suffices to verify it for any bounded continuous \(\phi_k(X), \, X\in {\cal X}\), due to the property of measures on $\mathbb R^d$ which are uniquely determined by the values of their integrals with continuous bounded functions (see, e.g., \cite[Theorem 1.2.4]{Krylov_ln}). 
The latter formula  (\ref{df22}) may be called one more  version of Dynkin's identity; it will be the basis for establishing existence. With a bit of abuse of the standard terminology,  \eqref{df22} may also be  called a martingale problem. Note, however, that weak uniqueness (= uniqueness in distribution) in this Theorem {\em given $(\mu_s, \, s\ge 0)$} does not mean a total uniqueness in distribution of the process under construction because there is no claim of uniqueness of \((\mu_s, s\ge 0)\), not even talking about a distribution in the space of trajectories. 

~

{\em Proof of Theorem \ref{thm1}}. For any \(n\ge 1\) consider a process \((X^{n}_t)\), with initial data \(X^{n}_0=X_0\) and intensities of jumps up and down, respectively, 
\[
\Lambda^+[t,X^{n}_{(t-1/n)_+}, \mu^{n}_{(t-1/n)_+}], \quad \Lambda^-[t,X^{n}_{(t-1/n)_+}, \mu^{n}_{(t-1/n)_+}].
\]
The process \((X^{n}_t)\) for each \(n\) is constructed by induction successfully on the  intervals \([0,1/n], [1/n, 2/n]\), etc. Due to the boundedness assumption on both intensities, there is no blow up and the processes for any \(n\) are defined for any \(t\ge 0\) as c\`adl\`ag processes without any point of jump accumulating. Moreover, for any \(t\) probability of jump exactly at time \(t\) for any \(X^{n}\) equals zero. 

Recall that the processes \((X^{n}_t,\, t\ge 0)\) for \(n\ge 1\) being constructed, we introduce on some probability space {\em independent}  equivalent processes  \((\xi^{n}_t, \, t\ge 0)\), and that  \(\mathbb E'\) stands in all cases  for the integration with respect to the {\em third variable}, e.g., 
\[
\mathbb E' \lambda^\pm(t,X_t^{n},\xi^{n}_t) := 
\int \lambda^\pm(t,X_t^{n},Y)\mu_t^{n}(dY).
\]
It can be checked that the assumptions of the Lemma \ref{lem1}   from the Appendix are satisfied. 

Indeed, given $\epsilon>0$, let us firstly choose $c_0>0$ so that 
\[
\mathbb P_{0,\mu_0}(|X_0|>c_0) < \epsilon/2. 
\]
On the event $(|X_0|\le c_0)$ we have for any $n$, 
\[
\sup_{0\le t\le T}(|x(X^n_t)| + |y(X^n_t)|) \le c_0+2T.  
\]
Further, since the intensity of jump up $\Lambda^+$ is bounded, say, 
$\Lambda^+ \le \bar\lambda$, then the number of jumps up on $[0,T]$ is bounded in probability, that is, uniformly with respect to $n$ (recall that $X^n_0 = X_0$), 
\[
\mathbb P(\sup_{0\le t\le T}k(X^n_t) - k(X_0)>c) \to 0, \quad c\to\infty. 
\]
Hence, the first condition (\ref{sko1}) of the Lemma \ref{lem1} for the family of processes $(X^n)$ follows. 

~

To check the second condition (\ref{sko2}), note that 
\begin{align*}
\mathbb P(|k(X^n_t) - k(X^n_s)|>0) \le \mathbb P(\mbox{at least one jump on $[s,t]$})
 \\\\
\le 1-\exp(-\bar\lambda |t-s|) 
\le \bar\lambda |t-s|.
\end{align*}
Next, for $\epsilon > 2h \ge 2|t-s|$, 
\begin{align*}
\mathbb P(|x(X^n_t) - x(X^n_s)|+|y(X^n_t) - y(X^n_s)|>\epsilon;\; \mbox{no jumps on $[s,t]$}) = 0.  
\end{align*}
So, the equality (\ref{sko2}) for the family of processes $(X^n)$ follows as required. 


~

Hence, on some new probability space there exist equivalent 
processes 
\((\tilde X^{n}_t, \tilde \xi^{n}_t)\), and a limiting pair \((\tilde X^{}_t, \tilde \xi^{}_t)\) such that for some subsequence \((\tilde X^{n'}_t, \tilde \xi^{n'}_t) \stackrel{\mbox{\small $\mathbb P$}}{\to}(\tilde X^{}_t, \tilde \xi^{}_t), \, n'\to \infty,\) for each \(t\). It follows due to the boundedness of all intensities that the limiting process \((\tilde X^{}_t, \tilde \xi^{}_t)\) is also stochastically continuous. More than that,  with probability one the pair $(\tilde X_t, \tilde \xi_t)$ is a process with a finite number of jumps on any bounded interval.  Moreover,  the property  \(\lim_{h\downarrow 0} \sup_n \sup_{t,s\le T; \, |t-s|\le h}\mathbb P(|\tilde X^n_t - \tilde X^n_s|>\epsilon) = 0\) implies that for any \(\epsilon>0\) there is a following convergence in probability, 
\[
\tilde X^{n'}_{(t-1/n')_+} \stackrel{\mbox{\small  $\mathbb P$}}{\to} \tilde X_t, \quad n'\to\infty. 
\]
For the sequel, denote by $\tilde {\cal F}_t^n$ the sigma-algebra $\sigma(\tilde X^n_s: \, 0\le s\le t)$, and again slightly abusing notations we will drop the upper index $n$ here.

The analogue of Dynkin's formula \eqref{df22} for the pair $(\tilde X_t^{n'}, \tilde \xi_t^{n'})$ reads, 
\begin{align}\label{df3}
&\displaystyle 
\mathbb E_{0,X_0} \left[\left(g(\tilde X^{n'}_{t_{m+1}}) -  g(\tilde X^{n'}_{t_m}) -   \int\limits_{t_m}^{t_{m+1}} \mathbb E' L(s,\tilde X^{n'}_{(s-1/n')_+},\tilde \xi^{n'}_{(s-1/n')_+})g(\tilde X^{n'}_s)\,ds\right)
 \right. \nonumber 
 \\ \\ 
 \nonumber 
&\displaystyle  
\left. \hspace{3cm} \times  \prod_{k=1}^{m} 
\phi_k(\tilde X^{n'}_{t_k})\right] = 0, \qquad t_1<\ldots < t_m < t_{m+1}.
\end{align}
The formula (\ref{df3}) follows straightforward from the ``complete expectation'' arguments (cf., for example, \cite{VZ}) and from the definition of intensities. 

Indeed, irrespectively on whether or not the intensities depend on the current state of the process ($\tilde X_t$), or on some past values with a delay, they intuitively mean that for $\delta>0$ we have,  
\begin{align*}
\!\mathbb P\left(\tilde X^n_{t+\delta}\! =\! (k+1,x',y+\delta), \!0\!\le \!x'\!\le\! \delta) | \tilde {\cal F}^n_t; \tilde X^n_{t}\! =\! (k,x,y), \tilde X^n_{(t-1/n)_+}\! =\! (k',x',y')\right) 
 \\\\
= \mathbb E' \Lambda^+(t,(k',x',y'),\tilde \xi^n_{(t-1/n)_+}) \delta + o(\delta), 
\end{align*}
and if $k(\tilde X^n_{(t-1/n)_+})>0$, 
\begin{align*}
\!\mathbb P\left(\tilde X^n_{t+\delta}\! =\! (k-1,x+\delta,y'), \, 0\!\le\! y'\!\le\! \delta | \tilde {\cal F}^n_t; \tilde X^n_{t}\! =\! (k,x,y), \tilde X^n_{(t-1/n)_+} \!= \!(k',x',y')\right)
 \\\\ 
= \mathbb E'\Lambda^-(t,(k',x',y'),\tilde \xi^n_{(t-1/n)_+})  \delta + o(\delta),
\end{align*}
and finally, 
\begin{align*}
\mathbb P\left(\tilde X^n_{t+\delta} = (k,x+\delta,y+\delta) | \tilde {\cal F}^n_t; \tilde X^n_{t} = (k,x,y),  \tilde X^n_{(t-1/n)_+} = (k',x',y')\right) 
 \\\\
=1- \mathbb E'(\Lambda^++ \Lambda^-)(t,(k',x',y'),\tilde \xi^n_{(t-1/n)_+})  \delta + o(\delta).
\end{align*}
Continuity of the intensities were implicitly used here; later this assumption will be dropped. 

~

Hence,  we can write for any bounded function $g$ in the domain of the operator $L$, 
\begin{align*}
\mathbb E \left(g(\tilde X^n_{t+\delta})|\tilde  {\cal F}^n_{t}\right)
= g(\tilde X^{n,+}_{t})\mathbb E'\Lambda^+(t,\tilde X^n_{(t-1/n)_+},\tilde \xi^n_{(t-1/n)_+}) \delta 
 \\\\ 
+g(\tilde X^{n,-}_{(t-1/n)_+})\mathbb E'\Lambda^-(t,\tilde X^n_{t},\tilde \xi^n_{(t-1/n)_+})\delta 
 \\\\
+ g(\tilde X^n_t+\delta)\left(1-\mathbb E'(\Lambda^++ \Lambda^-)(t,\tilde X^n_{(t-1/n)_+},\tilde \xi^n_{(t-1/n)_+}) \delta\right) + o(\delta), 
\end{align*}
as $\delta \downarrow 0$.
Therefore, 
\begin{align*}
\mathbb E \left(g(\tilde X^n_{t+\delta}) - g(\tilde X^n_{t})| \tilde {\cal F}^n_{t}\right)
 \\\\
= g(\tilde X^{n,+}_{t})\mathbb E'\Lambda^+(t,\tilde X^n_{(t-1/n)_+},\tilde \xi^n_{(t-1/n)_+})\delta 
 \\\\
+ 
g(\tilde X^{n,-}_{t})\mathbb E'\Lambda^-(t,\tilde X^n_{(t-1/n)_+},\tilde \xi^n_{(t-1/n)_+})\delta 
 \\\\
+ g(\tilde X^n_t+\delta)\left(1- \mathbb E'(\Lambda^++ \Lambda^-)(t,\tilde X^n_{(t-1/n)_+},\tilde \xi^n_{(t-1/n)_+}) \delta\right) - g(\tilde X^n_{t}) + o(\delta)
 \\\\
= \mathbb E'L(t,\tilde X^n_{(t-1/n)_+},\tilde \xi^n_{(t-1/n)_+})g(\tilde X^n_t)\delta + o(\delta).  
\end{align*}
	Applying now still intuitively the principle of ``complete probabilities'' (here more accurately it could be called ``complete expectation''), that is, taking a summation we can obtain the equation (\ref{df3}). Indeed, let us split the interval $[t_m, t_{m+1}]$ into $N$ equal small sub-intervals, $t_m = s_0 < \ldots < s_N = t_{m+1}$ so that $(s_{i+1}-s_i)=: \delta$  and write down,
\begin{align}
\mathbb E \left(g(\tilde X^n_{t_{m+1}}) - g(\tilde X^n_{t_{m}}) | \tilde{\cal F}^n_{t_m}\right)
= \sum_{i=0}^{N-1} \mathbb E \left(\mathbb E \left(g(\tilde X^n_{s_{i+1}}) - g(\tilde X^n_{s_i}) | \tilde {\cal F}^n_{s_i}\right)| \tilde {\cal F}^n_{t_m}\right) 
 \nonumber \\\nonumber\\\nonumber
=\sum_{i=0}^{N-1} \left[\mathbb E 
\left(\mathbb E'L(s_i,\tilde X^n_{(s_i-1/n)_+},\tilde \xi^n_{(s_i-1/n)_+})g(\tilde X^n_{s_i}) | | \tilde {\cal F}^n_{t_m}\right)\delta + o(\delta)\right] 
 \\\nonumber\\
= \mathbb E \left(\int_{t_m}^{t_{m+1}} \mathbb E'L(s,\tilde X^n_{(s-1/n)_+},\tilde \xi^n_{(s-1/n)_+})g(\tilde X^n_{s})\,ds |\tilde  {\cal F}^n_{t_m}\right) + o(1), \label{intuitive} 
\end{align}
which implies the equation (\ref{df3}) (recall that we dropped prime at $n$ to simplify notations). This intuitive calculus may be made strict; for the convenience of the reader we provide such a rigorous version in the Appendix 2. 

~

However, another easier way is just to recall the definition of intensity via the martingale property. In our particular case -- with jumps of the first component of $X$ just up or down -- the term intensity is applied to the random variable $\Lambda = \Lambda^\pm(t,\omega)$ if and only if for any bounded measurable function $g(X)$ ($X = (n,x,y)$)  with bounded derivatives with respect to $x$ and $y$ and for any $t_0\ge 0$, the process defined by the expression
\[
M^n_t:= g(\tilde X^n_t) - g(\tilde X^n_{t_0}) - \int_{t_0}^t \mathbb E'L(s,\tilde X^n_{(s-1/n)_+},\tilde \xi^n_{(s-1/n)_+})g(\tilde X^n_{s}) \, ds, \quad t\ge t_0,
\]
is a martingale: see, e.g., \cite[Sec.3.III.5.5]{LSh} for pure jump processes and for indicator functions, which extends straightforward to our case and to Borel measurable functions; the intuition behind this definition has been offered in the little calculus above (see also the Appendix 2). So, 
\[
\mathbb E\left [g(\tilde X^n_{t_{m+1}}) \!-\! g(\tilde X_{t_m}) \!-\! \!\!\int\limits_{t_m}^{t_{m+1}}\!\mathbb E'L(s,\tilde X^n_{(s-1/n)_+},\tilde \xi^n_{(s-1/n)_+})g(\tilde X^n_{s})\, ds | \tilde {\cal F}^n_{t_m}\right] \!= \!0 \; \mbox{(a.s.)}. 
\]
Therefore, it follows that 
\begin{align*}
\displaystyle 
\mathbb E_{0,X_0} \!\!\left(\!g(\tilde X^{n'}_{t_{m+1}})\! -\!  g(\tilde X^{n'}_{t_m})\! -\! \! \! \int\limits_{t_m}^{t_{m+1}} \mathbb E' L(s,\tilde X^{n'}_{(s-1/n')_+},\tilde \xi^{n'}_{(s-1/n')_+})g(\tilde X^{n'}_s)\,ds\!\right)\!
\!\prod_{k=1}^{m} 
\phi_k(\tilde X^{n'}_{t_k})
 \\\\
\displaystyle 
= \mathbb E_{0,X_0} \left(\prod_{k=1}^{m} 
\phi_k(\tilde X^{n'}_{t_k})\right) 
\mathbb E \left[g(\tilde X^n_{t_{m+1}})\right.
 \\ \\
\left. - g(\tilde X_{t_{m}}) - \int\limits_{t_{m}}^{t_{m+1}}\mathbb E' L(s,\tilde X^{n'}_{(s-1/n')_+},\tilde \xi^{n'}_{(s-1/n')_+})g(\tilde X^n_s)\, ds |\tilde  {\cal F}^n_{t_m}\right] = 0, 
\end{align*}
as required. This justifies the equation (\ref{df3}). 

Further, by continuity of \(\lambda\) and \(h\), and due to the stochastic continuity of the  processes \(\tilde X\) and \(\tilde \xi\), and since all integrand expressions are bounded, and by virtue of Lebesgue's bounded convergence Theorem, we obtain from \eqref{df3} in the limit with {\em continuous bounded} functions \((\phi_k)\), 
\begin{equation}\label{df44}
\mathbb E_{0,X_0} \left(g(\tilde X^{}_{t_{m+1}}) -  g(\tilde X^{}_{t_m}) -   \int\limits_{t_m}^{t_{m+1}} \mathbb E' L(s,\tilde X^{}_s,\tilde \xi^{}_s)g(\tilde X^{}_s)\,ds\right)
\prod_{k=1}^{m} \phi_k(\tilde X^{}_{t_k}) = 0. 
\end{equation}
Since the distribution of the random variable \(\tilde \xi_t\) is the same as the one of \(\tilde X_t\) -- let us denote it by \(\tilde \mu_t\) -- then \eqref{df44} can be equivalently written as  
\begin{equation}\label{df55}
\mathbb E_{0,X_0} \left(g(\tilde X^{}_{t_{m+1}}) -  g(\tilde X^{}_{t_m}) -   \int\limits_{t_m}^{t_{m+1}}  Lg(s,\tilde X^{}_s,\tilde \mu^{}_s)\,ds\right)
\prod_{k=1}^{m} \phi_k(\tilde X^{}_{t_k}) = 0. 
\end{equation}
As was mentioned earlier, due to the properties of measures on \(\mathbb R^d\) the formula \eqref{df55} holds true for any Borel bounded functions \((\phi_k)\), too. 
Due to \cite{Davis}, solution of the ``martingale problem'' \eqref{df55} -- or, more precisely, of the martingale problem 
\begin{equation}\label{mp7}
M_t:= g(\tilde X^{}_{t}) -  g(\tilde X^{}_{0}) -   \int\limits_{0}^{t}  Lg(s,\tilde X^{}_s,\tilde \mu^{}_s)\,ds, \; \; t\ge 0,  \quad \mbox{\em is a martingale}, 
\end{equation}
with a given family of marginal measures \((\tilde \mu_s, s\ge 0)\)  is unique. Hence,  according to \cite{Krylov-selection}, or  \cite[Theorem 4.4.2]{EthierKurtz} the limiting process \(\tilde X\) is Markov. The form of its generator 
with the required  intensities \(\Lambda^\pm\) follows from \eqref{df55}. This finishes the proof of the Theorem \ref{thm1}. 

\subsection{Weak uniqueness}\label{sec:2.2}
Emphasize that we will use essentially boundedness of all intensities and the condition that they are (uniformly) bounded away from zero. While it is clear that the boundedness from above may be  relaxed for the purpose of establishing existence -- e.g., under Lyapunov type conditions, or under a linear growth, or otherwise, -- and that boundedness away from zero is not required for the existence at all, yet for the uniqueness both boundedness from above and from below seems essential  (although also could be, apparently, slightly relaxed). On the other hand, continuity of the intensities in this section is not necessary and they are not assumed.

\begin{theorem}\label{thm2}
Let the assumptions (A1)--(A2) and (A4) be satisfied. Then, for any fixed distribution \({\cal L}(X_0)\),  there exists no more than one distribution of the process \((X_t, \, t\ge 0)\) with required intensities \(\Lambda^+[t,X,\mu_t]\) and   \(\Lambda^-[t,X,\mu_t]\). 
\end{theorem}
Recall that no Lipschitz assumptions on the intensities are assumed. In the calculus the total variation metric will be used. 

~

\noindent
Let 
\(\bar \Lambda[t,X,\mu] := \Lambda^+[t,X,\mu] + \Lambda^-[t,X,\mu]\).

~

{\em Proof of Theorem \ref{thm2}} 
is based on Skorokhod--Girsanov's change of measure formula for jump processes (see, e.g., \cite{LiptserShiryaev}). Suppose there are two solutions, \((X^1_t, \mu^1_t)\) and \((X^2_t, \mu^2_t)\). Denote by \(\Omega_n\) the event that the trajectory \(X\) has precisely \(n\) jumps on \([0,T]\).
Recall -- see, e.g., \cite{LiptserShiryaev}, \cite{Sko} -- that on the interval of time \([0,T]\) the density of one distribution with respect to the other -- we denote them by $\mathbb P^{\mu^i}, \, i=1,2$ -- on a typical trajectory $\omega = (t_1^\pm, \ldots, t_n^\pm)$ with overall $n\ge 0$ jumps up ($t_i^+$) or down ($t_j^-$) reads, 
\begin{eqnarray*}
&\displaystyle \rho_T:= \frac{d\mathbb P ^{\mu^2}}{d\mathbb P ^{\mu^1}}(\omega)|_{\Omega_n} 
 \\\\
&\displaystyle = \prod_{i=1}^{n}\frac{\Lambda^{\pm}[t^\pm_i,X_{t_i}, \mu^2_{t_i}]}{\Lambda^{\pm}[t_i^\pm,X_{t_i}, \mu^{1}_{t_i}]} \, \exp\left(-\int_0^T (\bar \Lambda^{}[t,X_{t}, \mu^2_{t}] - \bar \Lambda^{}[t,X_{t}, \mu^1_{t}])\,dt\right), 
\end{eqnarray*} 
where \(X = (X_s, \, 0\le s\le T)\) and \((t^\pm_i)\) are the moments of jumps of the trajectory \(X\), up or down, respectively; we keep the same sign at $\Lambda$, too, i.e., $\Lambda^+[t^+, \ldots]$ or, respectively, $\Lambda^-[t^-, \ldots]$. 
The usual convention \(\prod_{i=1}^{0} \ldots = 1\) is assumed.   Note that, of course, the number of jumps \(n\) is random -- i.e., it is a function of the trajectory -- but in any case it is almost surely finite due to the boundedness of the intensities. Note also that the expression $\rho_T$ above is a probability density. We have, 
\begin{eqnarray*}
&\displaystyle \mathbb E^{\mu^1} \prod_{i=1}^{n}\frac{\Lambda^{\pm}[t^\pm_i,X_{t_i}, \mu^2_{t_i}]}{\Lambda^{\pm}[t_i^\pm,X_{t_i}, \mu^1_{t_i}]} \, \exp\left(-\int_0^T (\bar\Lambda^{}[t,X_{t}, \mu^2_{t}] - \bar\Lambda^{}[t,X_{t}, \mu^1_{t}])\,dt\right)
 \\\\
&\displaystyle \!= \!\sum_{n=0}^{\infty}  \mathbb E^{\mu^1} 1(\Omega_n) \prod_{i=1}^{n}\frac{\Lambda^{\pm}[t_i^\pm,X_{t_i}, \mu^2_{t_i}]}{\Lambda^{\pm}[t_i^\pm,X_{t_i}, \mu^1_{t_i}]}  \exp\left(\!-\!\int_0^T\! (\bar\Lambda^{}[t,X_{t}, \mu^2_{t}] \!-\! \bar\Lambda^{}[t,X_{t}, \mu^1_{t}])\,dt\!\right)
 \\\\
&\displaystyle =  \sum_{n=0}^{\infty} \;\;
\mathbb E^{\mu^2} \idotsint\limits_{0 < t_1<\cdots <t_n<T}  \prod_{i=1}^{n}\Lambda^{\pm}[t_i^\pm,X_{t_i}, \mu^2_{t_i}] \, \exp\left(-\int_0^T \bar\Lambda^{}[t,X_{t}, \mu^2_{t}]\,dt\right)
\, \prod_{i=1}^{n}dt_i
 \\\\
&\displaystyle = \sum_{n=0}^{\infty}\mathbb P^{\mu^2}(\Omega_n) = 1. 
\end{eqnarray*}
Note that given the initial state $X_0$, the value without expectation $\mathbb E^{\mu^2}$ here equals, actually, 
\[
\sum_{n=0}^{\infty} \;\;
\idotsint\limits_{0 < t_1<\cdots <t_n<T}  \prod_{i=1}^{n}\Lambda^{\pm}[t_i^\pm,X_{t_i}, \mu^2_{t_i}] \, \exp\left(-\int_0^T \bar\Lambda^{\pm}[t,X_{t}, \mu^2_{t}]\,dt\right)
\, \prod_{i=1}^{n}dt_i, 
\]
which itself equals identically one, 
while expectation $\mathbb E^{\mu^2}$ relates to integration of each term over $X_0$ if it is random. 
It is worthwhile to recall that the rule of the evolution of the trajectory \(X\) between the moments of jumps \(t_i\) is deterministic and linear with rate one for the continuous components, and the discrete component does not change between any two consequent jumps. 

Now, we want to estimate the distance  {\em in total variation} between two probability measures in the space of trajectories, \(\mu^1_{[0,T]}\) and \(\mu^2_{[0,T]}\) and then to use the inequality that the distance between the marginals of any two measures does not exceed the distance of the measures themselves, 
\[
\varphi_T:= \|\mu^1_{T} - \mu^2_{T}\|_{TV} \le 
\|\mu^1_{[0,T]}-\mu^2_{[0,T]}\|_{TV}= 2-2\mathbb E^{\mu^1} \left(\rho_T\wedge 1\right)=:\psi_T. 
\]
Now, the idea is to estimate the right hand side in the last term via \(\varphi\) and, hence, to show that, at least, for small values of \(T>0\) this value equals zero. If this is realized, then the claim that \(\varphi_t=0\) for $t\le T$, $t\le 2T$, etc., and, eventually, for all \(t\ge 0\) would follow by induction. In fact, we will be able to estimate the right hand side via another expression with \(\psi_T\) itself. 
Note, by the way, that although normally marginal distributions of any process may not determine the distribution in the space of trajectories, in our case with intensities it is, of course, the case which follows from \cite{Davis}, as mentioned already in the proof of the Theorem \ref{thm1}. 

The first goal is to find a suitable lower bound for the value \(\mathbb E^{\mu^1} \left(\rho_T\wedge 1\right)\) from below. Let us split it as follows:
\[
\mathbb E^{\mu^1} \left(\rho_T\wedge 1\right) = \sum_{n=0}^{\infty}\mathbb E^{\mu^1} 1(\Omega_n)\left(\rho_T\wedge 1\right). 
\]
Further, we have for \(n = 0\), 
\begin{eqnarray*}
&\displaystyle 
\mathbb E^{\mu^1} 1(\Omega_0)\left(\rho_T\wedge 1\right)
\\\\
&\displaystyle =\mathbb E^{\mu^1}
1(\Omega_0) \exp\left(-\int_0^T (\bar\Lambda[t,X_{t}, \mu^2_{t}] - \bar\Lambda[t,X_{t}, \mu^1_{t}])\,dt\right)\wedge 1
 \\\\
&\displaystyle \ge \exp(-\int_0^T \|\lambda\|\|\mu^2_{t} - \mu^1_{t}\|_{TV}\,dt)\mathbb E^{\mu^1} 1(\Omega_0) 
 \\\\
&\displaystyle
\ge 
\exp(-\|\lambda\|\, T\psi_T)\mathbb E^{\mu^1} 1(\Omega_0). 
\end{eqnarray*}
All norms like $\|\lambda\|$ are sup-norms (except for the total variation norm, which is always shown explicitly). We used the fact that \(|\bar \Lambda[t,X,\mu]|\le  \|\lambda\|\), and that \[|\bar\Lambda[t,X_{t}, \mu^2_{t}] - \bar\Lambda[t,X_{t}, \mu^1_{t}]| \le \|\lambda\| |\mu^1_{[0,t]}-\mu^2_{[0,t]}|_{TV} \le \|\lambda\| |\mu^1_{[0,T]}-\mu^2_{[0,T]}|_{TV}, \; 0\le t\le T.\]
Similarly for \(n\ge 1\), with a notation $\tilde \Lambda^\pm[t^\pm, \ldots] := \ln \Lambda^\pm[t^\pm, \ldots]$,
\begin{eqnarray*}
&\displaystyle 
\mathbb E^{\mu^1} 1(\Omega_n)\left(\rho_T\wedge 1\right)
\\\\
&\displaystyle =\mathbb E^{\mu^1}
1(\Omega_n)\left\{\prod_{i=1}^{n}
\frac{\Lambda^\pm[t^\pm_i,X_{t_i}, \mu^2_{t_i}]}{\Lambda^{\pm}[t^\pm_i,X_{t_i}, \mu^1_{t_i}]} \times \right.
 \\\\
&\displaystyle \left. \times  \exp\left(-\int_0^T (\bar\Lambda[t,X_{t}, \mu^2_{t}] - \bar\Lambda[t,X_{t}, \mu^1_{t}])\,dt\right)\right\}\wedge 1
 \\\\
&\displaystyle \ge \mathbb E^{\mu^1}
1(\Omega_n)
\exp(-\sum_{i=1}^{n}|\tilde \Lambda[t^\pm_i, X_{t_i}, \mu^2_{t_i}] - \tilde\Lambda[t^\pm_i, X_{t_i}, \mu^1_{t_i}]|) 
 \\\\
&\displaystyle  \times 
\exp\left(-\int_0^T |\bar\Lambda[t,X_{t}, \mu^2_{t}] - \bar\Lambda[t,X_{t}, \mu^1_{t}]|\,dt\right).
\end{eqnarray*}
Minimum with $1$ here was dropped after all multipliers were estimated from below by the values less than one. Further, since the derivative of $\ln x$ is bounded on any interval $0<a\le x\le b$, say, by a constant $K$, we have with $a=\inf \lambda(\ldots)=:\underline \lambda$ and $b = \|\lambda\|$, 
\begin{align*}
|\tilde \Lambda^\pm[t^\pm_i, X_{t_i}, \mu^2_{t_i}] - \tilde\Lambda[t^\pm_i, X_{t_i}, \mu^1_{t_i}]| 
\le K |\Lambda^\pm[t^\pm_i, X_{t_i}, \mu^2_{t_i}] - \Lambda^\pm[t^\pm_i, X_{t_i}, \mu^1_{t_i}]|
 \\\\
\le K \|\Lambda\|\|\mu^2_{t_i} - \mu^2_{t_i}\|_{TV}.
\end{align*}
Hence, 
\begin{eqnarray*}
&\displaystyle 
\mathbb E^{\mu^1} 1(\Omega_n)\left(\rho_T\wedge 1\right)
 \\\\
&\displaystyle \ge \mathbb E^{\mu^1}
1(\Omega_n)
\exp(-\sum_{i=1}^{n}K\|\Lambda\|\, \|\mu^2_{t_i} - \mu^1_{t_i}\|_{TV})
 \\\\
&\displaystyle \times \exp(-\int_0^T \|\lambda\|\|\mu^2_{t} - \mu^1_{t}\|_{TV}\,dt).
\end{eqnarray*}
(Here by definition \(t_0=0\).) Here the  infimum \(\underline\lambda\)  is positive by the assumption.  Thus, using the bound \(1-\exp(-a) \le a\) and estimates 
\(\mathbb E^{\mu^1} 1(\Omega_0)\le \exp(-\underline{\lambda}T)\) and 
\(\displaystyle \mathbb E^{\mu^1} 1(\Omega_n)\le \frac{(\|\lambda\|T)^n}{n!}\exp(-\underline{\lambda}T)\), \(n\ge 1\), we get, 
\begin{eqnarray*}
&\displaystyle 
\frac12 \psi_T = 1 - \sum_{n}^{}\mathbb E^{\lambda^1} 1(\Omega_n)\left(\rho_T\wedge 1\right)
 \\\\
&\displaystyle \le 
(1-\exp(-\|\lambda\|\, T\psi_T)) \mathbb E^{\lambda^1} 1(\Omega_0)
 \\\\
&\displaystyle + \sum_{n=1}^{\infty}
\mathbb E^{\lambda^1} 
1(\Omega_n) 
\left(1-\exp\left(-K\sum_{i=1}^{n}\|\Lambda\|\, \|\mu^2_{t_i} - \mu^1_{t_i}\|_{TV}\right) \times \right.
 \\\\
&\displaystyle \left. \times \exp\left(-\int_0^T \|\lambda\| \|\mu^2_{t} - \mu^1_{t}\|_{TV}\,dt\right)\right)
 \\\\
&\displaystyle \le \|\lambda\| T \psi_T 
\mathbb E^{\mu^1} 1(\Omega_0) 
 \\\\
&\displaystyle 
+ \sum_{n=1}^{\infty}
\mathbb E^{\mu^1} 
1(\Omega_n) 
\left(1-\exp\left(- n K\|\Lambda\|\, \psi_T - T \|\lambda\| \psi_T\right)\right)
 \\\\
&\displaystyle \le \psi_T \exp(-\underline \lambda T) \left(\|\lambda\| T + \sum_{n=1}^{\infty}
(n K\|\Lambda\|\, + T \|\lambda\|) \, \frac{(\|\lambda\| T)^n}{n!}\right)
 \\\\
&\displaystyle = T\psi_T \exp(-\underline \lambda T) \, 
\left(\|\lambda\| + \sum_{n=0}^{\infty}
((n+1) K\|\Lambda\|\, + T \|\lambda\|) \, \frac{(\|\lambda\|)^{n+1} T^n}{(n+1)!}\right).
\end{eqnarray*}
The series in the right hand side here converges and does not exceed some constant, say, \(C>0\), if \(T\le 1\). Hence, overall, we obtain, 
\[
0\le  \frac12 \psi_T \le C T \psi_T, \quad T\le 1. 
\]
This implies that 
\[
\psi_T =0, \quad T < (2C)^{-1}\wedge 1,  
\]
and, therefore, also 
\[
\varphi_T =0, \quad T < (2C)^{-1}\wedge 1,  
\]
as required. In other words, we have shown that the two marginal measures \(\mu^1_t\) and \(\mu^2_t\) coincide for all \(t< (2C)^{-1}\wedge 1\). 

~

Further, note the constant \(C\) in this calculus does not depend on the initial distribution of the process. Hence, using the Markov property of the process and repeating the same arguments on \([T,2T]\),  \([2T,3T]\), etc., by induction we conclude that 
\[
\psi_t =0, \quad t\ge 0,   
\]
and, therefore, also 
\[
\varphi_t =0, \quad t\ge 0,  
\]
as required. So, the two measures \(\mu^1\) and \(\mu^2\)  on the space of trajectories are equal. The Theorem \ref{thm2} is proved.

\section*{Appendix 1}\label{ap1}
The following celebrated Lemma  is stated for the convenience of the reader.
\begin{lemma}[Skorokhod \mbox{\cite[Ch.1, \S 6]{Sko}}] \label{lem1}
Let $\xi^n_t$ ($t\ge 0$, $n=0,1,\ldots$) be some $d$-dimensional  stochastic processes defined on some probability space and let for any $T>0$, $\epsilon > 0$ the following hold true: 
\begin{eqnarray}
& \displaystyle \lim_{c\to\infty} \sup_n \sup_{t\le T}\mathbb P(|\xi^n_t|>c) = 0, \label{sko1}
 \\
& \displaystyle \lim_{h\downarrow 0} \sup_n \sup_{t,s\le T; \, |t-s|\le h}\mathbb P(|\xi^n_t - \xi^n_s|>\epsilon) = 0. \label{sko2}
\end{eqnarray}
Then  there exists a subsequence $n'\to\infty$ and a new probability can be constructed with processes $\tilde \xi^{n'}_t, \, t\ge 0$ and $\tilde \xi^{}_t, \, t\ge 0$, such that all finite-dimensional distributions of $\tilde \xi^{n'}_{\cdot}$ coincide with those of $\xi^{n'}_{\cdot}$ and such that  for any $\epsilon>0$ and all $t\ge 0$, 
\[
\mathbb P (|\tilde \xi^{n'}_t - \tilde \xi^{}_t|> \epsilon) \to 0, \quad n'\to\infty. 
\]

\end{lemma}

\section*{Appendix 2}\label{ap2}
Let us show how the intuitive calculus leading to (\ref{intuitive}) may be performed more rigorously. We have, 
\begin{align*}
\mathbb P\left(\tilde X^n_{t+\delta}\! =\! (k+1,x',y+\delta), \!0\!\le \!x'\!\le\! \delta) | \tilde {\cal F}^n_t; \tilde X^n_{t},  \tilde X^n_{(t-1/n)_+}\right) |_{\tilde X^n_{t}=(k,x,y)}
 \\\\
=\mathbb P\left(\tilde X^n_{t+\delta}\! =\! (k+1,x',y+\delta), \, 0\!\le \!x'\!\le\! \delta); \right. 
 \\\\
\left. \mbox{precisely one jump up on $[t,t+\delta]$} | \tilde {\cal F}^n_t; 
\tilde X^n_{t}, \tilde X^n_{(t-1/n)_+}\right) + o(\delta)
 \\\\
= \int_t^{t+\delta}
\mathbb E' \Lambda^+(s,\tilde X^n_{(s-1/n)_+},\tilde \xi^n_{(s-1/n)_+}) 
 \\\\
\times \exp\left(-\int_t^{t+\delta} \mathbb E' \bar\Lambda(r,\tilde X^n_{(r-1/n)_+},\tilde \xi^n_{(r-1/n)_+}) \,dr\right)\,ds + o(\delta)
 \\\\
= \int_t^{t+\delta}
\mathbb E' \Lambda^+(s,\tilde X^n_{(s-1/n)_+}, \tilde \xi^n_{(s-1/n)_+}) \,ds + o(\delta);  
\end{align*}
if $k(\tilde X^n_{(t-1/n)_+})>0$, 
\begin{align*}
\mathbb P\left(\tilde X^n_{t+\delta} = (k-1,x+\delta,y'), \, 0\le y'\le \delta | \tilde {\cal F}^n_t; \tilde X^n_{t}, \tilde X^n_{(t-1/n)_+}\right) |_{\tilde X^n_{t}=(k,x,y)}
 \\\\ 
=\mathbb P\left(\tilde X^n_{t+\delta} = (k-1,x+\delta,y'), \, 0\le y'\le \delta; \right. 
 \\\\ 
 \left. \mbox{precisely one jump down on $[t,t+\delta]$} | \tilde {\cal F}^n_t; \tilde X^n_{t}, \tilde X^n_{(t-1/n)_+}\right) + o(\delta)
 \\\\ 
= \int_t^{t+\delta}
\mathbb E' \Lambda^-(s,\tilde X^n_{(s-1/n)_+},\tilde \xi^n_{(s-1/n)_+})
 \\\\
\times  \exp\left(-\int_t^{t+\delta} \mathbb E' \bar\Lambda(r,\tilde X^n_{(r-1/n)_+},\tilde \xi^n_{(r-1/n)_+}) \,dr\right)\,ds + o(\delta)
 \\\\
= \int_t^{t+\delta}
\mathbb E' \Lambda^-(s,\tilde X^n_{(s-1/n)_+},\tilde \xi^n_{(s-1/n)_+}) \,ds + o(\delta);  
\end{align*}
and finally, 
\begin{align*}
\mathbb P\left(\tilde X^n_{t+\delta} = (k,x+\delta,y+\delta) | \tilde {\cal F}^n_t; \tilde X^n_{t},  \tilde X^n_{(t-1/n)_+}\right) |_{\tilde X^n_{t}=(k,x,y)}
 \\\\
=\mathbb P\left(\tilde X^n_{t+\delta} = (k,x+\delta,y+\delta); \mbox{no jumps on $[t,t+\delta]$} | \tilde {\cal F}^n_t; \tilde X^n_{t},  \tilde X^n_{(t-1/n)_+}\right) + o(\delta)
 \\\\
=\exp\left(-\int_t^{t+\delta} \mathbb E' \bar\Lambda(r,\tilde X^n_{(r-1/n)_+},\tilde \xi^n_{(r-1/n)_+}) \,dr\right) + o(\delta)
 \\\\
=1 -\int_t^{t+\delta} \mathbb E' \bar\Lambda(r,\tilde X^n_{(r-1/n)_+},\tilde \xi^n_{(r-1/n)_+}) \,dr+ o(\delta).
\end{align*}
Note that unlike in the intuitive calculus earlier, we did not use any regularity conditions on the intensities here and that  up to $o(\delta)$ the formulae above are all exact, not approximate as in the earlier intuitive version. 
Hence,  we can write for any bounded function $g$ in the domain of the operator $L$ (in particular, continuous in the second and third components of the state variable), 
\begin{align*}
\mathbb E \left(g(\tilde X^n_{t+\delta})|\tilde  {\cal F}^n_{t}\right)
= g(\tilde X^{n,+}_{t}) \int_t^{t+\delta}
\mathbb E' \Lambda^+(s,\tilde X^n_{(s-1/n)_+},\tilde \xi^n_{(s-1/n)_+}) \,ds 
 \\\\ 
+g(\tilde X^{n,-}_{(t-1/n)_+})\int_t^{t+\delta}
\mathbb E' \Lambda^-(s,\tilde X^n_{(s-1/n)_+},\tilde \xi^n_{(s-1/n)_+}) \,ds  
 \\\\
+ g(\tilde X^n_t+\delta)\left(1 -\int_t^{t+\delta} \mathbb E' \bar\Lambda(r,\tilde X^n_{(r-1/n)_+},\tilde \xi^n_{(r-1/n)_+}) \,dr\right) + o(\delta), 
\end{align*}
as $\delta \downarrow 0$.
Therefore, it follows rigorously that
\begin{align*}
\mathbb E \left(g(\tilde X^n_{t+\delta}) - g(\tilde X^n_{t})| \tilde {\cal F}^n_{t}\right)
 \\\\
= \int_t^{t+\delta}\mathbb E'L(s,\tilde X^n_{(s-1/n)_+},\tilde \xi^n_{(s-1/n)_+})g(\tilde X^n_s)\,ds + o(\delta).  
\end{align*}
So, we can rigorously apply  the principle of ``complete probabilities'' (or, better ``complete expectation'') taking a summation as earlier in the intuitive version: let us split the interval $[t_m, t_{m+1}]$ into $N$ equal small sub-intervals, $t_m = s_0 < \ldots < s_N = t_{m+1}$ so that $(s_{i+1}-s_i)=: \delta$; then we get similarly to (\ref{intuitive}),
\begin{align*}
\mathbb E \left(g(\tilde X^n_{t_{m+1}}) - g(\tilde X^n_{t_{m}}) | \tilde{\cal F}^n_{t_m}\right)
= \sum_{i=0}^{N-1} \mathbb E \left(\mathbb E \left(g(\tilde X^n_{s_{i+1}}) - g(\tilde X^n_{s_i}) | \tilde {\cal F}^n_{s_i}\right)| \tilde {\cal F}^n_{t_m}\right)
 \\\\
= \mathbb E \left(\int_{t_m}^{t_{m+1}} \mathbb E'L(s,\tilde X^n_{(s-1/n)_+},\tilde \xi^n_{(s-1/n)_+})g(\tilde X^n_{s})\,ds |\tilde  {\cal F}^n_{t_m}\right) + o(1),  
\end{align*}
which rigorously implies the equation (\ref{df3}).

\section*{Acknowledgements}
The techniques used in this paper were stimulated by the methods developed in a long-term joint work on formally quite different McKean-Vlasov SDE equations with Yu. Mishura, as well as in fruitful discussions of the author on the same subject with D. \v{S}iska, and L. Szpruch. S. Pirogov, A. Rybko, and G. Zverkina  helped to find some (quite a few) technicalities to be corrected in the earlier versions of the text. The author is sincerely thankful to all these colleagues and to two referees for very useful advice. The deepest gratitude is to Professor Alexander Dmitrievich Solovyev (06.09.1927 -- 06.04.2001) who was the author's supervisor at BSc and MSc programmes at Moscow State University. 



\end{document}